\documentclass[12pt]{article}
\usepackage{latexsym}
\usepackage{amssymb}
\def\mymedskip{\vskip\medskipamount}
\def\mymedbreak{\par \ifdim\lastskip<\medskipamount
  \removelastskip \penalty-100 \mymedskip \fi}
\def\myaftermedspace{\par \ifdim\lastskip<\medskipamount
  \removelastskip \penalty55\mymedskip\fi}
\newcommand{\eop}{{\unskip\nobreak\hfil\penalty50
          \hskip2em\hbox{}\nobreak\hfil$\Box$
          \parfillskip=0pt \finalhyphendemerits=0 \par}}
\newenvironment{proof}%
{\mymedbreak{\noindent\bf Proof:\enspace}}{\eop\myaftermedspace}
{\mymedbreak{\noindent\bf Proof of Theorem #1:\enspace}}{\eop\myaftermedspace}
\newenvironment{example}%
{\mymedbreak\refstepcounter{Exc}
                      {\em Example \theExc:}\enspace}%
{\eop\myaftermedspace}
\newtheorem{teor}{Theorem}[section]
\newcounter{Exc}

\newtheorem{defi}[teor]{Definition}
\newtheorem{lem}[teor]{Lemma}
\newtheorem{cor}[teor]{Corollary}
\newtheorem{con}[teor]{Conjecture}
\newtheorem{prop}[teor]{Proposition}
\newtheorem{rem}[teor]{Remark}
\newcommand{\beq}{\begin{equation}}
\newcommand{\eeq}{\end{equation}}
\newcommand{\beql}[1]{\begin{equation} \label{#1}}
\newcommand{\eeql}{\end{equation}}
\newcommand{\beqa}{\begin{eqnarray*}}
\newcommand{\eeqa}{\end{eqnarray*}}
\newcommand{\beqal}[1]{\begin{eqnarray} \label{#1}}
\newcommand{\eeqal}{\end{eqnarray}}
\newcommand{\beqan}{\begin{eqnarray}}
\newcommand{\eeqan}{\end{eqnarray}}
\newcommand{\bpf}{\begin{proof}}
\newcommand{\epf}{\end{proof}}

\newcommand{\bex}[1]{\begin{example} \label{#1}}
\newcommand{\eex}{\end{example}}

\newcommand{\cI}{{\cal I}}
\newcommand{\cL}{{\cal L}}

\newcommand{\bF}{{\bf F}}

\newcommand{\per}{{\rm per}}
\newcommand{\chr}{{\rm char}}

\newcommand{\gs}{\sigma}

\newcommand{\gw}{\omega}

\newcommand{\ga}{\alpha}
\newcommand{\gb}{\beta}

\newcommand{\gl}{\lambda}

\newcommand{\ord}{{\rm ord}}
\newcommand{\lcm}{{\rm lcm}}

\newcommand{\bbF}{\mathbb{F}}
\newcommand{\bbZ}{\mathbb{Z}}
\newcommand{\bbE}{\mathbb{E}}

\newtheorem{fact}[teor]{Fact}
\newtheorem{problem}{Problem}
\newtheorem{exercise}{Exercise}

\newtheorem{alg}[teor]{Algorithm}
\newcommand{\ben}{\begin{enumerate}}
\newcommand{\een}{\end{enumerate}}
\newcommand{\bit}{\begin{itemize}}
\newcommand{\eit}{\end{itemize}}
\newcommand{\Tm}[1]{Theorem~\protect\ref{#1}}
\newcommand{\Le}[1]{Lemma~\protect\ref{#1}}
\newcommand{\Pn}[1]{Proposition~\protect\ref{#1}}

\newcommand{\btm}[1]{\begin{teor} \label{#1}}
\newcommand{\etm}{\end{teor}}
\newcommand{\btmn}[2]{\begin{teor}[#1] \label{#2}}
\newcommand{\etmn}{\end{teor}}
\newcommand{\ble}[1]{\begin{lem} \label{#1}}
\newcommand{\ele}{\end{lem}}
\newcommand{\bLe}[1]{\begin{Lemma} \label{#1}}
\newcommand{\eLe}{\end{Lemma}}
\newcommand{\bpn}[1]{\begin{prop} \label{#1}}
\newcommand{\epn}{\end{prop}}
\newcommand{\bde}[1]{\begin{defi} \label{#1}}
\newcommand{\ede}{\end{defi}}
\newcommand{\bco}[1]{\begin{cor} \label{#1}}
\newcommand{\eco}{\end{cor}}
\newcommand{\bcorn}[2]{\begin{cor}[#1] \label{#1}}
\newcommand{\ecorn}{\end{cor}}
\newcommand{\bcon}[1]{\begin{con} \label{#1}}
\newcommand{\econ}{\end{con}}
\newcommand{\bfact}[1]{\begin{fact} \label{#1}}
\newcommand{\efact}{\end{fact}}
\newcommand{\bpr}[1]{\begin{problem} \label{#1}}
\newcommand{\epr}{\end{problem}}
\newcommand{\bprnn}[1]{\begin{problemnn} \label{#1}}
\newcommand{\eprnn}{\end{problemnn}}
\newcommand{\bprn}[2]{\begin{problem}[#1] \label{#2}}
\newcommand{\eprn}{\end{problem}}
\newcommand{\bexer}[1]{\begin{exercise} \label{#1}}
\newcommand{\eexer}{\end{exercise}}
\newcommand{\bre}[1]{\begin{rem} \label{#1}}
\newcommand{\ere}{\end{rem}}
\newcommand{\balg}[1]{\begin{alg} \label{#1}}
\newcommand{\ealg}{\end{alg}}
\newcommand{\bqu}{\begin{question}}
\newcommand{\equ}{\end{question}}
\newcommand{\bs}{\begin{solution}}
\newcommand{\es}{\end{solution}}
\newcommand{\bh}{\begin{hint}}
\newcommand{\eh}{\end{hint}}
\newcommand{\bms}[1]{\begin{multisolution}{#1}}
\newcommand{\ems}{\end{multisolution}}
\newcommand{\bquo}{\begin{quote}}
\newcommand{\equo}{\end{quote}}
\newcommand{\la}{\langle}
\newcommand{\ra}{\rangle}

\newcommand{\et}{\tilde{e}}
\newcommand{\st}{\tilde{s}}
\newcommand{\tti}{\tilde{t}}
\providecommand{\keywords}[1]
{
  \small	
  \textbf{\textit{Keywords---}} #1
}
\begin{document}
\begin{titlepage}
\title{
Some basic results on finite linear recurring sequence subgroups} 
\date{December 1, 2020}
\author{Henk D. L.\ Hollmann\\University of Tartu\\Institute of Computer Science\\ Tartu
50409, Estonia.\\Email: {\tt henk.d.l.hollmann@ut.ee}
\and
Medet Zhanbulatuly\\
10 Anson Road,\\
 International Plaza, \#38-16, \\
 079903 Singapore\\
Email: {\tt medet.ntu@gmail.com}
}
\maketitle
\begin{abstract}
An {\em $f$-subgroup\/} is a linear recurring sequence subgroup, a multiplicative subgroup of a field whose elements can be generated (without repetition) by a linear recurrence relation, with characteristic polynomial~$f$. It is called {\em non-standard\/} if it can be generated in a  non-cyclic way (that is, not in the order $\ga^i, \ga^{i+1}, \ga^{i+2} \ldots$ for a zero $\ga$ of~$f$), and {\em standard\/} otherwise. 
We will show that a finite $f$-subgroup is necessarily generated by a subset of the zeros of $f$. We use  this result to improve on a recent theorem of Brison and Nogueira.  A old question by Brison and Nogueira asks if there exist {\em automatically non-standard\/} $f$-subgroups, $f$-subgroups that cannot be generated by a zero of~$f$. We answer that question affirmatively by constructing infinitely many examples. 
\end{abstract}
%
\end{titlepage}
%
%
\section{Introduction}
A finite multiplicative subgroup of $\bbF^*$, the nonzero elements in a commutative field $\bbF$, is necessarily cyclic (see, e.g.\ \cite[Chapter 1, Lemma 1]{weil}), that is, of the form $\la \ga \ra=\{1, \ga, \ga^2, \ldots\}$ for some field element~$\ga$ of finite order. 
In a sequence of papers \cite{bri, bn1, bn-mat,bn2, bn3, bn4, bn5}, Brison and Nogueira have investigated presentations of such a subgroup in a finite field by means of a (periodic) {\em linear recurring sequence\/}. 
Here, we say that a  linear recurring sequence $s$ {\em presents\/} (or {\em generates\/}) a finite subgroup~$M$ if $s$ has period $|M|$ and $M=\{s_0, s_1, \ldots, s_{|M|-1}\}$; moreover, if the linear recurrence relation satisfied by the sequence~$s$ has characteristic polynomial~$f$, then we refer to~$M$ as 
an {\em $f$-subgroup\/} and to the generating sequence~$s$  as an {\em $f$-sequence\/}.

For example, consider the (Fibonacci-type) linear recurrence relation 
\beql{LEFR} s_n=s_{n-1}+s_{n-2}\eeql
with characteristic polynomial $f(x)=x^2-x-1$. If $\ga$ is a zero of $f$ in some finite field~$\bbF_q$ of size $q$ and characteristic $p>0$ (so with $q$ a power of the prime $p$), then the sequence $s$ with $s_n=\ga^n$ for $n\in \bbZ$ obviously satisfies the recurrence relation (\ref{LEFR}), so if~$\ga$ has order $m$,
then $M=\{1,\ga, \ldots, \ga^{m-1}\}$ is a presentation of $M$ by $s$, and hence $M$ is an $f$-subgroup. Usually, this is essentially the only way to present $M$, and then we refer to $M$ as a {\em standard\/} $f$-subgroup.  However, it sometimes happens that there is another, essentially {\em different\/} way to present $M$. When $p\equiv \pm2\bmod 5$, the polynomial~$f$ is irreducible and $q=p^2$; by combining results from  \cite{bri} and \cite{bn1}, it can be shown that the subgroup $M\leq \bbF_q^*$ generated by a zero of~$f$ is a standard $f$-subgroup except when $p=3, q=9$: in that case, if $\gw\in\bbF_9\setminus \bbF_3$,  the sequence 
\[s=\ldots, 1, \gw, 1+\gw, 1-\gw, -1, -\gw, -1-\gw, -1+\gw, 1, \gw, \ldots\]
satisfies the recursion (\ref{LEFR}) and has period 8, so this $f$-sequence~$s$ presents $M=\bbF_9^*$. This gives 6 presentations for $M$ but $f$ has at most two zeros, and so 4 of these presentations are {\em non-standard\/}. We now refer to $M=\bF_9^*$ as a {\em non-standard\/} $f$-subgroup. (In the remaining cases where $p=5$ or $p\equiv \pm1\bmod 5$, both zeros $a, b$ of $x^2-x-1$ are in $\bbF_p$, so $q=p$; later in this paper 
we show that $M=\la a,b\ra$ is in fact the only possible non-standard $f$-subgroup, 
but we do not know if this possibility ever occurs.)

As the above example suggests, non-standard $f$-subgroups tend to be rare (at least when we put some restrictions on $f$, see Section~\ref{dis}). In the case where $f(x)\in \bbF_q[x]$ is {\em irreducible\/}, it is known that an $f$-subgroup is necessarily generated by a zero of~$f$. The cases where $f$ is irreducible of degree 2 have been fully classified, first when $q$ is prime in \cite{bn1}, then for general $q$ in \cite{cyclic-subm}, using a result from \cite{bn3}. Interestingly, it turns out that for irreducible polynomials $f$, non-standardness occurs precisely when the irreducible cyclic code related to $f$ has extra automorphisms, for more details we again refer to \cite{cyclic-subm}. Various types of non-standard irreducible cases have been identified, and we may  
well know them all. 


In the general case, a clear understanding of the relation between the zeros of a polynomial~$f$ and an $f$-subgroup was still missing. 
This relation is investigated in Section~\ref{gen}. The main result in this section states that an $f$-subgroup 
presented by a periodic sequence~$s$ is equal to the group generated by the zeros of the minimal polynomial~$f_s$ of~$s$, hence is generated by some subset of the zeros of~$f$.
Note that this generalizes the result for irreducible polynomials mentioned earlier. We use our results to remove the degree condition ``$k\leq p$'' in~\cite[Theorem 2.3]{bn5}.

Then in Section~\ref{auto}, we define a class of non-standard $f$-subgroups that we 
named {\em automatically non-standard\/}. These are subgroups $M$ of~$\bbF^*$ presented by a periodic $f$-sequence~$s$ 
where no zero of $f$ generates $M$. 
Note that if $s_{n+1}/s_n=\ga $ for all $n$, then 
$f(\ga)=0$ would follow from the  recurrence relation; so the non-standardness 
of such an $f$-subgroup is inherent from the definition. This class is interesting since, as we show in this section,  it turns out to be nonempty! We describe infinitely many examples, and we also prove an extension result.
In \cite[end of Observation 1.5]{bn1}, the authors state " \ldots if $f$ is irreducible [with zero $\xi$] then any $f$-subgroup has the form $\la \xi\ra$ (considered as
group), but we have no proof that this must occur in general." So our examples answer the implicit question here in the negative.
In \cite{bn0}, automatically non-standard $f$-subgroups in the complex field where named {\em non-standard of the second type\/}. The authors informed us that they also knew examples in finite fields already in 2015.

To make the paper self-contained, in Section~\ref{pre} we establish some notation and we quickly sketch a proof of most of the basic results concerning linear recurrence relations that are needed here. We basically follow the elegant approach in \cite{fm}, but some of our proofs may be new and our method may be of independent interest.

Finally, in Section~\ref{dis} we suggest some directions for further research and we discuss some open problems.

Part of this work is based on ideas obtained during a visit of the first author to Brison and Nogueira in 2013 and on later work by the second author \cite{meu,mef}.

\section{\label{pre}Preliminaries}
In this section, we quickly sketch the required background on linear recurrence relations and linear recurring sequences. To
describe the results, we will basically use the framework from \cite{fm}. While all the results in this section are known, parts of our approach may be new. 
%
For an element $\ga\neq 0$ in a field~$\bbF$, we write $\la \ga \ra$ to denote the group $\{1, \ga, \ga^2, \ldots\}$ generated by~$\ga$.
We will write $K\leq G$ to denote that $K$ is a subgroup of a group~$G$.
For other references for this material, see, for example, \cite{lak}, \cite{ln}, \cite{mce}, \cite{zier}.

Let $\bbF$ be an arbitrary field. A (two-way infinite) sequence $s=\{s_n\}_{n\in \bbZ}$ with  $s_n\in \bbF$ for all $n\in \bbZ$
is called a ($k$th-order)  {\em linear recurring sequence
over $\bbF$\/} if there exist  $c_0, \ldots, c_{k-1}$ in~$\bbF$, 
where $c_0\neq 0$ if $k>0$,
such that
\beql{LErr}s_n=c_{k-1}s_{n-1}+c_{k-2}s_{n-2}+\cdots +c_1s_{n-k-1}+c_0s_{n-k}\eeql
for all $n\in \bbZ$. A relation of the form (\ref{LErr}) is called a ($k$th-order)
{\em linear recurrence relation\/}. The polynomial 
\[f(x)=x^k-c_{k-1}x^{k-1}+\cdots-c_1x-c_0\]
in $\bbF[x]$ is referred to as the {\em characteristic polynomial\/} of (\ref{LErr}). We will say that the sequence~$s$ {\em satisfies the recursion $f$\/}, or that $s$ is an {\em $f$-sequence\/}, if (\ref{LErr}) holds for all $n\in \bbZ$, and we will write $\cL_{\bbF}(f)$ to denote the collection of all $f$-sequences $s$ over $\bbF$. We say that a sequence satisfies a recursion $f$ where $f(x)=cx^k+\cdots$ has leading coefficient $c\neq 0$ if $s$ satisfies $c^{-1}f$.
A sequence $s$ is called {\em cyclic} if there exists $\ga\in\bbF$ such that $s_{n+1}/s_n=\ga$ for all $n\in \bbZ$. Note that such a sequence is an $f$-sequence if and only if $f(\ga)=0$. 

\bre{LRext}\rm
Requiring instead that the coefficients $c_0, \ldots, c_{k-1}$ in (\ref{LErr}) are contained in some {\em extension field\/}~$\bbE$ of~$\bbF$ would not have widened the notion of a linear recurring sequence over~$\bbF$. That is, if a sequence $s$ with $s_n\in \bbF$ for all $n\in \bbZ$ satisfies {\em some\/} linear recurrence relation, then it also satisfies one where all coefficients are in~$\bbF$.  For periodic sequences, this follows from~\Tm{LTfs} below. In general, this follows from the existence of Berlekamp-Massey-type algorithms to compute the minimal polynomial, see, e.g., \cite{nor95,nor,sag}.
\ere

The collection $\cL_\bbF$ of all sequences over $\bbF$ forms a vector space under point-wise addition and scalar multiplication, defined by  $(s+t)_n=s_n+t_n$ and $(\gl s)_n=\gl s_n$  for $n\in \bbZ$ and $\gl\in \bbF$. 
The {\em (left) shift operator\/} $\gs$ operates on this vector space as $\gs(s)=t$, where $t_n=s_{n+1}$ for $n\in \bbZ$. Let $\gs^k$ denote $\gs$ composed with itself $k$ times,
and define $(a\gs^k+b\gs^l)(s)=a\gs^k(s)+b\gs^l(s)$ for $a,b\in \bbF$.  
The next proposition collects some basic results for this set-up. 
\bpn{LPmains}\rm \mbox{}\\
(i) 
The set of all operators $f(\gs)$ for $f(x)\in \bbF[x]$, under addition and composition, is isomorphic to the polynomial ring~$\bbF[\gs]$, and the above defines a ring-action of $\bbF[\gs]$ on the vector space~$\cL_\bbF$ (in technical terms, $\cL_\bbF$ is a left $\bbF[\gs]$-module).\\
(ii) 
A sequence $s\in \cL_\bbF$ 
is an $f$-sequence if and only if $f(\gs)s=0$.\\
(iii) 
If $s$ is an $f$-sequence and $t$ is a $g$-sequence, then $s+t$ is a $h$-sequence whenever both~$f,g$ divide $h$.\\
(iv) 
If $f(x)$ has degree $k$, then the collection $\cL_{\bbF}(f)$ of $s\in\cL_\bbF$ for which $f(\gs)s=0$ is a $k$-dimensional subspace of~$\cL_\bbF$.\\
(v)  For every $s\in \cL_\bbF$, the collection~$\cI_\bbF(s)$  
of polynomials $f(x)\in\bbF[x]$ for which $f(\gs)s=0$ is an {\em ideal\/} in~$\bbF[x]$, and so there is a unique monic polynomial $f_s(x)$ such that 
$f(\gs)s=0$ if and only if $f_s(x)$ divides~$f(x)$.\\
(vi) If $(f,g)=1$, then $\cL_\bbF(f)\cap \cL_\bbF(g)=\cL_\bbF(1)=\{0\}$.
\epn
\bpf (Sketch) (i) -- (iii) are easy;  for (i), remark that indeed $f(\gs)(g(\gs)s)=(f(\gs)g(\gs)) s$ for all~$f(x), g(x)\in\bbF[x]$.  Part (iv) follows from the observation that a sequence $s\in\cL_\bbF$ satisfying a polynomial $f(\gs)$ of degree $k$ is completely determined by $(s_0, \ldots, s_{k-1})$. The polynomial $f_s(\gs)$ in (v) is easily seen to be the monic polynomial with smallest degree in~$\cI_\bbF(s)$.
Finally, to see (vi), note that if $(f,g)=1$ in $\bbF[x]$, then there are $a(x), b(x)\in \bbF[x]$ for which $a(x)f(x)+b(x)g(x)=1$. Now if $s$ is both an $f$-sequence and a $g$-sequence,  then by (ii) $s$ is also a $1$-sequence, hence $s=0$.
\epf 
We will refer to the polynomial $f_s$ in (v) above as the {\em minimal polynomial\/} or the {\em minimal recursion\/} for $s$. In view of (i), in what follows we will identify $\bbF[x]$ and $\bbF[\gs]$, and we will use $x$ instead of $\gs$ to denote the left shift operator.

Next, we will describe the general solution of a linear recurrence relation.
To this end, as in \cite{rio} we define the binomial coefficients ${n\choose j}$ for $n,j\in \bbZ$ 
by the relations
${n \choose j} = {n -1\choose j}+{n-1 \choose j-1}$
and 
$ {n+j-1\choose j}=(-1)^j {-n\choose j}$
for $n, j= 1, 2, 3, \ldots$ and in addition,
${n\choose 0}=1$ for $n\in\bbZ$
and ${0 \choose j} =0$
for $j\in \bbZ$, $j\neq 0$.
It follows that
${n \choose -j}=0$ for $n\in \bbZ,   j=1,2,\ldots$ and
${n\choose n+j}=0$ for $n=0,1,2,\ldots, j=1, 2,\ldots$.
We sketch a proof of the following.
%
\btm{LTmainrecrel} \rm Suppose that $f(x)\in \bbF[x]$ has 
degree $m$ and factors completely over an extension $\bbE$ of~$\bbF$. If $f$ has  $t$ distinct zeros $\ga_1, \ldots, \ga_t\neq 0$ in $\bbE$, where $\ga_i$ has multiplicity~$e_i$ ($1\leq i \leq t$), then the collection $\cL_\bbE(f)$ of $f$-sequences over~$\bbE$ consists of the sequences  $s=\{s_n\}_{n\in \bbZ}$ for which
\[s_n=\sum_{i=1}^t\ga_i^n\sum_{j=0}^{e_i-1} c_{i,j}{n \choose j}\qquad (n\in \bbZ)\]
for suitable $c_{i,j}\in\bbE$.
\etm
\bpf
By \Pn{LPmains}, part (iv), the vector space $\cL_\bbE(f)$ of $f$-sequences over $\bbE$ has dimension $m=e_1+\cdots +e_t$, so we need to find $m$ independent solutions. As a consequence of \Pn{LPmains}, part (vi), it is sufficient to prove the theorem for $f(x)=(x-\ga)^m$ with $\ga\neq 0$, which follows if we show that the $m$ sequences $s^{(j)}$ with $s_n^{(0)}=1, s^{(1)}_n={n\choose 1}, \ldots, s^{(m-1)}_n={n\choose m-1}$ for $n\in \bbZ$ constitute $m$ independent solutions for the recursion $f(x)=(x-1)^m$.

Since $(s^{(j)}_0, \ldots, s^{(j)}_{m-1})$
has the first nonzero entry equal to a 1, in position $j$, it is evident that the sequences  $s^{(j)}$ for $j=0, \ldots, m-1$ are independent.
To complete the proof, we have to show that each sequence $s^{(j)}$ with $0\leq j<m$ satisfies the recurrence relation
$s_n - {m\choose 1}s_{n-1}+\cdots + (-1)^{m-1}{m\choose m-1} s_{n-m+1}+(-1)^m s_{n-m}=0$, that is,
we must show that $S(n,j)=\sum_{i=0}^m (-1)^{i} {m\choose i} {n-i \choose j}=0$
for $j=0, \ldots, m-1$ and for $n\geq m$. 

In fact, 
\[S(n,j)=
\left\{
\begin{array}{ll}
0, & \mbox{if $0\leq j\leq m-1$};\\
{n-m\choose j-m},& \mbox{if $j\geq m$}.
\end{array},
\right.
\]
Indeed, this is \cite[Chapter 1, (5a)]{rio}, or it follows from 
\[\sum_{j\geq0}S(n,j)x^j=(1+x)^{n}(1-1/(x+1))^m=(1+x)^{n-m}x^m=
\sum_{j\geq m}{n-m\choose j-m}x^{j}.\]
\epf
%


\bre{LRcharp}\rm
The advantage of describing the solutions in terms of binomial coefficient sequences ${n\choose j}$ instead of the conventional $n^j$  is that the binomial solutions also work over fields with finite characteristic $p$ when $e_i> p$, cf.\ \cite[Remark 8.23]{ln}. 
\ere

For completeness' sake, we also mention the following result.
\btmn{\mbox{\protect\cite[Theorem 8.27]{ln}}}{LTfseqper}Let $f(x)\in \bbF[x]$ with $f(0)\neq 0$, where $\chr(\bbF)=p>0$. Then every $f$-sequence~$s$ satisfies $\per(s)|\ord(f)$.
\etmn

Finally, we need an expression for the minimal recursion $f_s$ of a periodic sequence $s$.
Let  $u^{(m)}$ be the periodic sequence with period~$m$ for which $(u^{(m)}_0, \ldots, u^{(m)}_{m-1})=(0,\ldots, 0, 1)$. First we show the following.
\ble{LLures}\rm The sequence $u^{(m)}$ has minimal recursion $x^m-1$.
\ele 
\bpf
Since  $(x^m-1) u^{(m)}=0$, the minimal recursion divides $x^m-1$. 
Conversely, 
from~(2) we immediately see that any 
recursion for~$u^{(m)}$ must  have order~$k\geq m$.
\epf
\ble{LLsperm}\rm The sequence $s$ is periodic with period~$m$ if and only if $s=\st(x)u^{(m)}$ 
with $\st(x)=s_0x^{m-1}+\cdots+s_{m-2}x+s_{m-1}$.
\ele
\bpf
Obvious from the definition of~$u^{(m)}$.
\epf
The next theorem can also be derived using the approach in \cite[Theorem 8.25]{ln}.
\btm{LTfs}\rm
Let $s$ be periodic with period $m$ over $\bbF$. 
Then 
\beql{LEfs}f_s(x)=(x^m-1)/(x^m-1, \st(x)),\eeql
where $\st(x)$ is defined as in \Le{LLsperm}; in particular, $f_s(x)\in \bbF[x]$.
\etm
\bpf
By  \Le{LLsperm} we have that $s=\st(x)u^{(m)}$; then using \Le{LLures}, we conclude that $f(x)s=0$ if and only if 
$x^m-1$ divides $f(x) \st(x)$. The latter condition holds if and only if every zero of $x^m-1$ with multiplicity $e$ that occurs with multiplicity $h<e$ in $\st(x)$ occurs in $f(x)$ with multiplicity at least $e-h$. This is the case precisely when 
the right hand side of~(\ref{LEfs}) divides $f(x)$; now (\ref{LEfs}) follows from \Pn{LPmains}, part (v). Finally, note that the polynomial $f_s(x)$ defined by (\ref{LEfs}) is automatically contained in~$\bbF[x]$. 
\epf
%
The following basic result on periodic sequences can
be obtained by combining various known results, but we prefer to give a simple direct proof.
%
\btm{LTfundfseq} \rm 
%
Let $f(x)=\prod_{i=1}^t (x-\ga_i)^{e_i}$ be a polynomial in a field $\bbF$ with $p=\chr(\bbF)$, with distinct zeros $\ga_1, \ldots, \ga_t\neq0$ in some extension $\bbE$ of $\bbF$, and let $s=\{s_n\}_{n\in \bbZ}$ be a nonzero $f$-sequence of 
period $m$ in $\bbF$, for which $(m,p)=1$ if $p>0$.
Define $J=\{i\in\{1, \ldots, t\}\mid \ga_i^m=1\}$. Then the following hold.\\
(i) The sequence $s$ has minimal recursion $f_s$, where $f_s$ is as defined in~\Tm{LTfs}. In particular, $f_s$ divides $f$ 
and $f_s$ has no multiple zeros.\\
(ii) There are $c_{i}\in \bbE$   ($i\in J$) such that 
\beql{LEsgen}s_n=\sum_{i\in J}c_{i} \ga_i^n \qquad (n\in \bbZ);\eeql
moreover, $f_s(x)=\prod_{\{i\mid c_i\neq 0\}}(x-\ga_i)$.\\
(iii) $m=\lcm(\ord(\ga_i) \mid i\in J\}=\lcm(\ord(\ga_i) \mid \mbox{\rm $\ga_i$  is zero of $f_s$} \}$. 
\etm
\bpf 
(i) Under the conditions on $m$, the polynomial $x^m-1$ has no multiple zeros, hence the claim is a direct consequence of \Pn{LPmains}, part (v) and \Tm{LTfs}.
%

\noindent
(ii) Since $f_s$ divides $f$ and $f_s(x)|x^m-1$ (see part (i)), we conclude that 
every zero of~$f_s$ is of the form $\ga_i$ with $i\in J$. Since $s$ is an $f_s$-sequence, the expression (\ref{LEsgen})  follows from \Tm{LTmainrecrel}. Finally, again by \Tm{LTmainrecrel}, a sequence $s$ of the form  (\ref{LEsgen}) satisfies the recurrence $\prod_{\{i\mid c_i\neq 0\}}(x-\ga_i)$ but not one of smaller degree.
%

\noindent
(iii) Let  $r=\lcm(\ord(\ga_i)\mid i\in J)$ and  $r'=\lcm(\ord(\ga_i) \mid \mbox{\rm $\ga_i$   is zero of $f_s$} \}$. By part (ii), by~(\ref{LEsgen}), and by the definition of $f_s$, every $\ga_i$ actually occurring in the expression (\ref{LEsgen}) for~$s_n$ is in fact a zero of $f_s$, hence $r'$ is a period of $s$, and as a consequence $m|r'$. Conversely, since $\ga_i^m=1$ for every $i\in J$, we have $\ord(\ga_i)|m$ for $i\in J$, hence $r'|r|m$. We conclude that $r=m$. 
\epf
%
%
\section{\label{gen}General results for $f$-subgroups}
We now investigate the relation between the characteristic polynomial $f$ of a recurrence relation and the period of an $f$-sequence presenting an $f$-subgroup. The following fundamental result may be considered as a generalization of~\cite[Lemma 2.1]{bn2} and~\cite[Lemma 2.2]{bn5}.  
\btm{LCfundfseq}\rm Let $f(x)\in \bbF[x]$ and let $M$ be a finite $f$-subgroup of size~$m$ in some extension of~$\bbF$, presented by an $f$-sequence $s$ with minimal recursion $f_s$.  Then 
$f_s$ has no multiple zeros and $M$
is the group generated by the zeros of~$f_s$, and hence is generated by a subset of the zeros of~$f$.
\etm
\bpf
If $s$ presents $M$, where $M$ has size $m$, then $s$ is periodic with smallest period $m$, and $(|M|,\chr(\bbF))=1$ if $\chr(\bbF)>0$. Now the claims follow 
from~\Tm{LTfundfseq}.
\epf
\bex{LEdeg2}\rm Here we show how the above can be used to investigate $f$-subgroups where $f\in\bbF[x]$ has degree 2. Suppose that $f$ has two distinct zeros $a,b\in \bbF$ (for the case $a=b$, see \Tm{LTbn5} below).
By \Tm{LCfundfseq},
an $f$-subgroup is one of 
$\la a,b\ra$, $\la a\ra$, or~$\la b\ra$. Moreover, if the $f$-subgroup is smaller than $\la a, b\ra$, then by the same theorem, any $f$-sequence presenting it must have a minimal recursion of degree 1, so is cyclic; hence  the $f$-subgroup is standard. We conclude that $\la a,b\ra$ is the only candidate to be a non-standard $f$-subgroup.
%
\eex
\bre{LRperbin} \rm Let $p$ be a prime. 
Suppose that $h>\nu_p(j!)$, where $ \nu_p(n)$ denotes the largest power of $p$ dviding $n$.  Then for every integer $u$ we have that 
\beql{LEfac}{p^h u +n\choose j}=\frac{(p^hu+n)(p^hu+n-1)\cdots(p^hu+n-j+1)}{j!}\equiv 
{n\choose j}\bmod p.\eeql
As a consequence, if every zero $\ga_i$ of the minimal polynomial $f_s(x)$ of a sequence $s$ satisfies $\ga_i^m=1$, then $(m,p)=1$ and by \Tm{LTmainrecrel} the sequence $s$ has a period $p^hm$ (this need not be the minimal period). This fact can be used for an alternative approach to \Tm{LTfundfseq}, along the lines of \cite[Theorem 3.2.1.]{mef}.
Note that (\ref{LEfac}) and \Tm{LTfundfseq} together show that a linear recurring sequence is periodic in characteristic $p$ precisely when either $p>0$ or when $p=0$ and every zero of its minimal recursion has finite order and multiplicity 1. 
%

%
For completeness' sake, we remark that a much more precise result is available. 
Indeed, \cite[Lemma 1]{tre} (see also \cite[Section 4]{fray}) implies that if $p^{r-1}\leq j<p^r$, then the sequence $s^{(j)}$ in~$\bbF_p$ with $s^{(j)}_n={n\choose j}$ for $n\in\bbZ$ has minimal period~$p^r$. As a consequence, we can obtain an alternative proof for~\Tm{LTfseqper}. To this end, note that if $q=p^s$ and $f(x)\in \bbF_q[x]$ has $\per(f)=n=n_0p^r$ with $(n_0,p)=1$, then the general expression for an~$f$-sequence~$s$ in~\Tm{LTmainrecrel} immediately implies that such a sequence has $n_0p^r=n=\ord(f)$ as a period: indeed, $f(x)|x^n-1=(x^{n_0}-1)^{p^r}$, hence $n_0$ is a period of all the $\ga_i$, and $p^r$ is a period of all the binomial coefficients that occur. 
\ere

%

%
%
We now present an application of the above results.
In \cite[Theorem 2.3]{bn5}, the authors showed that an $f$ subgroup for $f(x)=(x-a)^k\in \bbF_{p^e}$ is standard provided that $k\leq p$. We will use the above results to eliminate the extra condition on~$k$ and to simplify the proof.
\btm{LTbn5}\rm Let $a\in \bbF_q^*$,  where $q=p^e$ with $p$ a prime, and let $f(x)=(x-a)^k\in \bbF_q[x]$ with~$k\geq1$. Then an $f$-subgroup $M$ is necessarily of the form $M=\la a\ra$ and is standard as an $f$-subgroup.
\etm
\bpf
Suppose that $s$ is an $f$-sequence
with smallest period $m$ such that $M=\{s_0=1, \ldots, s_{m-1}\}$ has size $m$ and is an $f$-subgroup in some extension $\bbE$ of $\bbF_q$. Then we may assume that $|\bbE|$ is finite, and since $m$ divides $|\bbE|-1$, we have $(m,p)=1$. The minimal polynomial $f_s(x)$ of $s$ divides $f(x)$, hence is of the form $(x-a)^e$, and by  \Tm{LTfundfseq}, we have $e=1$, hence 
$s$ is cyclic with $s_{n+1}/s_n=a$ ($n\in \bbZ$) and $M=\la a\ra$ is standard.
\epf
\section{\label{auto}Automatically non-standard $f$-subgroups}
%
In this section, we answer a question that was implicitly raised in \cite{bn1}, Observation 1.5, namely whether it is always true that an $f$-subgroup is generated by a zero of~$f$. We need some preparation.

For positive integers $m$, define the {\em $m$th cyclotomic polynomial\/} $\phi_m(x)$ inductively by letting 
\beql{LEcycl} x^m-1=\prod_{d|m} \phi_d(x).\eeql
So, for example, $\phi_1(x)=x-1$, $\phi_2(x)=x+1$, $\phi_3(x)=x^2+x+1$, and $\phi_4(x)=x^2+1$. 
It is well-known and not difficult to prove(see, e.g.,  \cite{gar}) that in fact every $\phi_m(x)$ is a polynomial with integer coefficients, that is, $\phi_m(x)\in \bbZ[x]$, and $\deg(\phi_m)=\varphi(m)$, where $\varphi$ is the Euler function defined by $\varphi(m)=|\bbZ^*_m|$, the number of integers $k$ with $1\leq k<m$ for which $(k,m)=1$.

Now let $\bbF$ be a field with $\chr(\bbF)=p$, and let $m$ be a positive integer with $(m,p)=1$ if $p>0$, so that $x^m-1$ has no multiple zeros. The {\em $m$th cyclotomic polynomial\/} $\phi_{p,m}$ over~$\bbF$ is the polynomial $\phi_m$, reduced modulo~$p$ if $p>0$. By its definition, $\phi_{p,m}$ has as its zeros precisely those $m$th roots of unity that have order~$m$, that is, the {\em primitive $m$th roots of unity\/}. Indeed, let $\ga$ be a zero of $\phi_{p,m}$ in a suitable extension of~$\bbF$. By assumption,  $(m,p)=1$, so $\ga$ has order~$m$.
Then the zeros of $x^m-1$ are $\ga^k$ ($0\leq k<m$), and $\ga^k$ is a primitive root of unity precisely when $(k,m)=1$. As a consequence, 
\[\phi_{p,m}(x)=\prod_{k\in \bbZ_m^*}(x-\ga^k),\] 
where the product is over all integers $k$ with $1\leq k<m$ for which $(k,m)=1$. 

Now assume that $f(x)\in\bbF[x]$, where $p=\chr(\bbF)$, and that $M$ is an $f$-subgroup of (finite) size $m$, that is, there exists an $f$-sequence $s$ of minimal period $m$ in some finite extension $\bbE$ of $\bbF$ such that $M=\{s_0, \ldots, s_{m-1}\}$. 
Then $(m,p)=1$ if $p>0$, and  \Tm{LTfundfseq} applies: the sequence $s$ has  a 
minimal polynomial $f_s(x)$ in $\bbF[x]$ that has 
no multiple zeros, where $f_s$ divides $f$ and $M$ is generated by the zeros of $f_s$. 
Now if $s$ is cyclic, say $s_{n+1}/s_n=\ga$ ($n\in \bbZ$), then $M=\la \ga \ra$, so $\ga$ has order $m$, and since~$s$ satisfies $f_s$ we see from (\ref{LErr}) that $f_s(\ga)=0$, and hence $f(\ga)=0$. But this cannot happen if $f$ does not have a primitive $m$th root of unity 
as a zero. Observe also that since $M=\{s_0, \ldots, s_{m-1}\}$ is a group, we have that $s_0+\cdots +s_{m-1}=0$, hence $x-1|\st(x)$ (with $\st(x)$ as defined in \Le{LLsperm}). With~\Tm{LTfs} in mind, this motivates the following definition. 
%
%
\bde{LDans}\rm Let $\bbF$ be a field with $p=\chr(\bbF)$, and let 
$m>1$ be an integer satisfying 
$(m,p)=1$ if $p>0$.\\
(i) If the polynomial $f(x)$ in~$\bbF[x]$ divides $(x^m-1)/(x-1)\phi_{p,m}(x)$ 
and if  $s$ is an $f$-sequence for which $M=\{s_0, s_1, \ldots, s_{m-1}\}$ is a subgroup of size $m$ in some extension of~$\bbF$,
then we refer to 
$M$ as an {\em automatically non-standard $f$-subgroup\/}.\\
\noindent (ii) If $s_0, \ldots, s_{m-1}$ are such that $M=\{s_0, \ldots, s_{m-1}\}$ is a subgroup of size~$m$ in some extension of~$\bbF$ and if $\phi_{p,m}(x)$ divides $\st(x)=s_{m-1}+s_{m-2}x+\cdots +s_0x^{m-1}$, then we say that the subgroup $M$ is {\em automatically non-standard\/}.
\ede
By the discussion preceding this definition, the following should come as no surprise. 
\btm{LTans} (i) An automatically non-standard subgroup $M=\{s_0, \ldots, s_{m-1}\}$ is automatically non-standard $f$-subgroup with $f(x)=(x^m-1)/(x^m-1,\st(x))$. \\
(ii) An automatically non-standard $f$-subgroup is a non-standard $f$-subgroup.
\etm
\bpf
(i) The periodic sequence $s$ with period $m$ defined by $s_0, \ldots, s_{m-1}$ generates $M$ and has minimal polynomial $f_s(x)=f(x)$ by~\Tm{LTfs}. By definition, $\phi_{p,m}|\st(x)$, and since $M$ is a group also $s_0+\cdots+s_{m-1}=0$ and $x-1|\st(x)$; since $m>1$, we conclude that  $f(x)|(x^m-1)/(x-1)\phi_{p,m}(x)$ so $M$ is automatically non-standard $f$-subgroup.\\
(ii) 
Let the automatically non-standard $f$-subgroup~$M$ have order~$m>1$. By the definition of $\phi_{p,m}$, 
no zero of $(x^m-1)/(x-1)\phi_{p,m}(x)$, hence no zero of $f$, generates~$M$. 
Hence $M$ cannot be generated by a cyclic $f$-sequence, so must be non-standard. (See \cite[Lemma 1.3 (b)]{bn4} for more details.)
\epf
%
%
A priori, it is not evident that automatically non-standard objects even exist. Moreover, we have the following simple  negative result.
\btm{LTnspp}\rm A subgroup $M$ with $|M|=r^e$ with $r$ prime is never automatically non-standard. 
\etm
\bpf
Let $M$ be an $f$-subgroup, for some polynomial $f$. By \Tm{LTfundfseq}, we may assume that $f$ has distinct zeros $\ga_1, \ldots, \ga_k$ and that $M=\la\ga_1, \ldots, \ga_k\ra$. Every $\ga_i$ has order~$r^h$ for some integer~$h\leq e$, hence $M$ can only be generated by the~$\ga_i$ if some~$\ga^i$  has order $r^e$, that is, if it generates $M$. 
\epf
As a consequence, there are no automatically non-standard subgroups of sizes $2,3,4,5, 7, 8, 9, 11, 13, 16$ and $17$, but there could be automatically non-standard subgroups of size 6, 10, 12, 14, and 15.
%
%
%
%
\bex{LEns6}\rm Suppose that $M$ is an automatically non-standard multiplicative subgroup of size 6 in some finite field $\bbF_q$ of characteristic $p$ (so  with $p\neq 2, 3$). Let $M=\{s_0=1, \ldots, s_5\}$, where $\st(x)=x^5+s_1x^4+\cdots +s_5$ is a multiple of 
the polynomial $\phi_6(x)=x^2-x+1$ containing all primitive $6$th roots of unity. With $\ga^2=\ga-1$ we have
 $\ga^3=-1, \ga^4=-\ga, \ga^5=-\ga+1$, and so $M=\{1,\ga, \ga-1, -1, -\ga, -\ga+1\}$;
moreover, $\st(x)\equiv0\bmod x^2-x+1$ and $\st(x)\equiv 0\bmod x-1$ (since $M$ is a subgroup), hence 
\[(-x+1)+s_1(-x)+s_2(-1)+s_3(x-1)+s_4x+s_5=0,\]
or equivalently, 
\[1-s_2-s_3+s_5=0, \qquad  -1-s_1+s_3+s_4=0,\]
and 
\[1+s_1+s_2+s_3+s_4+s_5=0.\]
By a careful examination of all possibilities, it can be shown \cite{mef} that necessarily the characteristic $p$ of the field satisfies $p=7$, with $M=\bbF_7^*$ and $s=(1, 3, 4, 6, 5, 2)$ (or one of the other 5 sequences such as $(1,5,3,4,2,6)$ obtained from this one by multiplying by $s_i^{-1}$ and shifting); the first sequence $s$ has $\st(x)=x^5+3x^4+4x^3+6x^2+5x+2=(x^3-2x^2+2x-1)(x^2+5x+5)$ and $(x^2+5x+5, x^6-1)=1$; hence this $s$ is an $f$-sequence with $f(x)=(x+1)(x^2+x+1)=x^3+2x^2+2x+1$
and generates the automatically non-standard $f$-subgroup~$\bbF_7^*$.
%
\eex
\noindent The above example can be generalized as follows.
\btm{LTans}\rm Let $m=p-1=2r$ with $r\geq3$ odd and $p$ prime (that is, $p$ prime, $p\equiv 3\bmod 4$ with $p\geq 7$). Then the subgroup $M=\bbF_p^*$ is automatically non-standard with respect  to the polynomial $f(x)=(x+1)(x^r-1)/(x-1)$.
\etm
\bpf
Let $m=p-1=2r$ with $r\geq3$ odd, and let $s$ be the sequence with period $m=p-1=2r$ defined by
\[s_0, \ldots, s_{m-1}=1,-2, 3, -4, \ldots, -(r-1), r, -r, r-1, -(r-2), \ldots, 2,-1.\]
Then 
\[ \st(x)=x^{2r-1}-2x^{2r-2}+\cdots -(r-1)x^{r+1}+rx^{r} -rx^{r-1}+(r-1)x^{r-2} +\cdots +2x-1.\]
Let 
\[a(x)=1-2x+3x^2-\cdots -(r-1)x^{r-2}+rx^{r-1}.\]
Then $a(x)=b'(x)$ with 
\[b(x)=x-x^2+x^3-\cdots -x^{r-1}+x^r=x(x^r+1)/(x+1)=(x^{r+1}+x)/(x+1),\]
hence 
\[a(x)=((r+1)x^r+1)(x+1)-(x^{r+1}+x)\cdot 1 )/(x+1)^2=(rx^{r+1}+(r+1)x^r+1)/(x+1)^2.\]
and 
\beqa \st(x)&=&-a(x)+x^{2r-1}a(x^{-1})\\
&=&((-rx^{r+1}-(r+1)x^r-1) +(rx^r+(r+1)x^{r+1}+x^{2r+1}))/(x+1)^2\\
&=&(x^{2r+1}+x^{r+1}-x^r-1)/(x+1)^2=(x^{r+1}-1)(x^r+1)/(x+1)^2.
\eeqa
%
Now if $\xi$ is a primitive $m$th root of unity in $\bbF_p$, then $\xi^{2r}=1$ but $\xi^{r}\neq 1$ (since $r\geq3$) and $\xi^2\neq 1$, so $\xi^r=-1$ and $\xi$ is a zero of $(x^r+1)/(x+1)$. As a consequence, $\st(x)$ contains all primitive $m$th roots of unity and $M=\{s_0, \ldots, s_{m-1}\}=\bbF_p^*$ is automatically non-standard in $\bbF_p$.

Since $r$ is odd, we have that
\beqa (x^{2r}-1, \st(x))&=&(x^{2r}-1,  (x^{r+1}-1)(x^r+1)/(x+1)^2) \\
&=&((x+1)(x^r-1), (x^{r+1}-1)/(x+1))(x^r+1)/(x+1)\\
&=&(x-1)(x^r+1)/(x+1),
\eeqa
hence by \Tm{LTfs}, the minimal polynomial for the sequence $s$ is 
\[f_s(x)=(x^{2r}-1)/(x^{2r}-1, \st(x))=
(x+1)(x^r-1)/(x-1).\]
\epf
\bex{LEmed11}\rm
There is another way to show that for $m=p-1=2r$ with $r\geq3$ odd the group $\bbF_p^*$ is automatically non-standard.
Indeed, let $s$ be the sequence with period $m=p-1$ defined by
\[s_0, \ldots, s_{m-1}=1,-1, 3, -3, \ldots, r-2, -(r-2); r; 2, -2, 4, -4, \ldots, r-1,-(r-1); r+1.\]
Then 
\[ \st(x)=x^{2r-1}-2x^{2r-2}+\cdots -(r-1)x^{r+1}+rx^{r} -rx^{r-1}+(r-1)x^{r-2} +\cdots +2x-1.\]
In \cite{mef} it is shown that  in fact 
\beqa \st(x)&=&[(x+1)(x^{r-2}+3x^{r-4}+5x^{r-6}+\cdots+(r-2)x) +r](x-1)(x^r+1)/(x+1)\\
&=&a(x)(x-1)(x^r+1)/(x+1)\eeqa
with $a(x)=(x+1)(x^{r-2}+3x^{r-4}+5x^{r-6}+\cdots+(r-2)x) +r$
 (it is not too difficult to check this directly).
Note that a primitive $(2r)$th root of unity $\xi$ satisfies $\xi^r=-1$ and $\xi\neq -1$ for $r\geq 3$, hence every generator of $M=\bbF_p^*$  is a zero of $\st(x)$; as a consequence, $\bbF_p^*$ is automatically non-standard. 
In order to compute the minimal polynomial $f_s(x)$ for the sequence $s$, we use \Tm{LTfs}. 
We have that
\[(x^{2r}-1, \st(x))=((x-1)(x^r+1)/(x+1)) ((x+1)(x^r-1)/(x-1),a(x)).\]
Now $(x+1,a(x))=(x+1,r)=1$, and 
it is not difficult to verify \cite{mef} that 
\[(x^2-1)a(x)-(x^2+1)(x^r-1)/(x-1)=-(r-1)x-(r+1)=-(r-1)(x-1/3).\]
Hence
\[((x+1)(x^r-1)/(x-1),a(x)) =(x^r-1, x-1/3).\]
Now $1/3$ is zero of $x^r-1$ when $3^r=3^{(p-1)/2}=1$ in $\bbF_p$, that is, when $3$ is a (nonzero) square in $\bbF_p$.
If $p\neq 2,3$, then by quadratic reciprocity, see, e.g., \cite{ir}, this happens precisely when $p\equiv \pm 1 \bmod 12$. 
In our case, $p\equiv 3 \bmod 4$ and $p\geq 7$, so $p\not\equiv 1\bmod 12$ (and certainly not $p\equiv 3 \bmod 12$) and we have that
\[f_s(x)=\left\{
\begin{array}{ll}
(x+1)(x^r-1)/(x-1),&\mbox{if $p\equiv  7\bmod 12$};\\ 
(x+1)(x^r-1)/(x-1)(x-1/3),&\mbox{if $p\equiv 11\bmod 12$}.
\end{array}
\right.
\]
\eex
\section{\label{ext}Extension for automatically non-standard subgroups}
In this section, we will prove the following {\em extension\/} theory.
\btm{LTANSE}\rm
If $M$ is an automatically non-standard $f$-subgroup in some field~$\bbE$,
then every finite multiplicative subgroup of~$\bbE$
that contains~$M$ is again automatically non-standard, with respect to the polynomial $g(x)=f(x^k)$.
\etm
\bpf
By our assumptions on~$M$, there exists an $f$-sequence~$s$  with period $m=|M|$ over $\bbE$ such that $M=\{s_0, \ldots, s_{m-1}\}$, where $f(x)|(x^m-1)/((x-1)\phi_{p,m}(x)$. Since the minimal recursion $f_s$ for $s$ satisfies $f_s(x)|f(x)$, we conclude from \Tm{LTfs} that $(x-1)\phi_{p,m}(x)|\st(x)$, where $\st(x)$ is the polynomial associated with the sequence~$s$.
%
%
Now let $L$ be a subgroup with $M\leq L$
and $|L|=k|M|$, say. Suppose that $e_0, e_1, \ldots, e_{k-1}$ are a system of distinct coset representatives of $L/M$. Now consider the presentation for $L$ as 
\[e_0s_0, e_1s_0, \ldots, e_{k-1}s_0, e_0s_1, \ldots, e_0s_{m-1}, \ldots, e_{k-1}s_{m-1},\]
that is, $L=\{t_0, \ldots, t_{km-1}\}$ with $t_{kj+i}=e_is_j$ for $0\leq i<k, 0\leq j<m$. For the associated polynomial $\tti(x)$ of the sequence $t$ (extended with period $km$) we now find
\beqa \tti(x)&=&\sum_{u=0}^{km-1}t_ux^{km-1-u}
= \sum_{i=0}^{k-1}\sum_{j=0}^{m-1}e_is_j x^{mk-1-jk-i}\\
&=&  \sum_{i=0}^{k-1} e_ix^{k-1-i}\sum_{j=0}^{m-1}s_jx^{k(m-1-j)}
=\et(x)\st(x^k),
\eeqa
where $\et(x) =\sum_{i=0}^{k-1} e_ix^{k-1-i}$.
Now since $\st(x^k)|\tti(x)$ and $\phi_{p,km}(x)|\phi_{p,m}(x^k)|\st(x^k)$,  we conclude from~\Tm{LTfs} that $f_t(x)|f_s(x^k)$, hence $f_t(x)|f(x^k)$. Finally, using again that $\phi_{p,km}(x)|\phi_{p,m}(x^k)$, we have that $f(x^k)|(x^{km}-1)/(x^k-1)\phi_{p,m}(x^k)|(x^{km}-1)/(x-1)\phi_{p,km}$. So with $g(x)=f(x^k)$, the subgroup $L$ is presented by the $g$-sequence~$t$, hence $L$ is an automatically non-standard $g$-subgroup.
\epf
\bre{LRlimit}\rm
So far, we have only constructed automatically non-standard groups $\bbF_p^*$ of size $p-1$ for primes $p\equiv 3\bmod 4$, and their extensions (in the sense of \Tm{LTANSE}) of size $m=k(p-1)$ with $k>0$ an integer with $(k,p)=1$. Note also that this type of extension to an automatically non-standard group $M=\bbF_q^*$ for $q$ a prime power is not always possible, since for example $M=\bbF_9^*$, of size $8=2^3$, is ruled out by \Tm{LTnspp}.
%
%
However, there seems to be no obvious reason why there cannot be an automatically non-standard subgroup of every size~$m$ of the form $m=rs$ with $r,s>1$ and $(r,s)=1$, although we do not know an example in every such case. Taking $M=\bbF_{11}^*$ provides a automatically non-standard group  for $m=10$.  Extension of the $m=6$ example in characteristic $p=7$ by taking $k=2$ gives a non-standard group of size $12$ in $\bbF_{7^2}^*$, so the smallest undecided case is $m=14$. This size may still be small enough to be handled by an exhaustive search, if needed with the help of a computer.
\ere
\section{\label{dis}Discussion and open problems}
%
%
Trivially, there are no non-standard subgroups of order $m\leq 3$; however {\em every\/} subgroup $M\leq \bbF^*$ of size $m\geq4$ is a non-standard $f$-subgroup for $f(x)=(x^m-1)/(x-1)=x^{m-1}+\cdots +x+1$ in {\em any\/} characteristic $p$, provided that $(m,p)=1$. Indeed, let $M=\la \ga \ra$ with $\ga$ in some extension $\bbF_q$ of $\bbF_p$ (such an element exists if $(m,p)=1$), and let $\pi$ be a permutation of $1, 2, \ldots, m-1$. Then $s_0=1, s_1=\ga^{\pi(1)}, \ldots, s_{m-1}=\ga^{\pi(m-1)}$ is a presentation of $M$ by a periodic $f$-sequence $s$ with period $m$. Indeed, if $s_n, s_{n+1}, \ldots, s_{n+m-2}$ are $m-1$ distinct elements of $M$, then $s_{n+m-1}=-s_n-s_{n+1}-\cdots -s_{n+m-2}$ is the remaining element in $M$ distinct from $s_n, \ldots, s_{n+m-2}$.
So there are at least $(m-1)!$ $f$-sequences $s$ with $s_0=1$ presenting $M$, and there are only $\varphi(m)$ presentations $1,\gb, \gb^2,\ldots$ with $\gb=\ga^k$, $(k,m)=1$, a primitive $m$th root of 1 in $M$. Now $(m-1)!>m-1\geq \varphi(m)$ for $m\geq 4$. For example, if $m=4$ and $p$ is odd, there exists a primitive 4th root of unity, say~$\ga$,  in~$\bbF_p^*$ (if $p\equiv 1\bmod 4$) or in~$\bbF_{p^2}\setminus \bbF_p$ (if $p\equiv 3\bmod 4$). By \Tm{LTfs}, both presentations $s_0, s_1, s_2,s_3=1, -1, \ga, -\ga$ and $1, \ga, -\ga, -1$ have minimal polynomial $f_s(x)=x^3+x^2+x+1=(x+1)(x^2+1)$.

In order to avoid such slightly trivial examples, we need to put further constrants on the polynomial $f$. One possibility is to require that $f(x)$ is irreducible over some field $\bbF_q$, and to take $M=\la \ga \ra$ for some zero $\ga$ of $f(x)$ in view of \Tm{LTfundfseq}.
Another possibility is to require that no zero of $f(x)$ generates $M$, which leads to the automatically non-standard subgroups considered here in Sections~\ref{auto} and~\ref{ext}.

In view of the above, and given \Tm{LTfundfseq}, the most ambitious goal would be to determine, for every prime $p$, all {\em minimally non-standard polynomials\/} in characteristic~$p$, that is, all polynomials $f$ over a field of characteristic~$p$
for which the subgroup $M$ generated by the zeros of $f$ can be presented as $M=\{s_0, s_1, \ldots, s_{m-1}\}$ with $m=|M|$ for
a non-cyclic sequence $s$ with minimal period $m$ and minimal recursion $f_s=f$ (an even more ambitious goal would be to count the {\em number\/} of such presentations). 
Note that every minimally non-standard polynomial has degree at least 2, as every polynomial $f(x)=x-a$ of degree 1 is of course standard. 

With this point of view, all minimally non-standard polynomials $f(x)\in\bbF_q$ of degree~2 and of the form $f(x)=(x-\xi)(x-\xi^q)$ for some $\xi\in\bbF_{q^2}\setminus \bbF_q$ (that is, irreducible over $\bbF_q$) have been determined in~\cite{cyclic-subm}; there is a relation with irreducible cyclic codes having extra automorphisms. In \cite{bn4} and in \cite{cyclic-subm}, various other classes of irreducible non-standard polynomials are described; obviously these are all minimally non-standard. 
So besides classifying all irreducible non-standard polynomials (which might be possible at least for degree at most 3), now one of the main open problems is the determination of  the minimally non-standard polynomials $f(x)\in\bbF_q$ of degree~2, of the form $f(x)=(x-a)(x-b)$ with $a, b\in\bbF_q\setminus\{0,1\}$, $a\neq b$, with  $M=\la a,b\ra$. To the best of our knowledge,  the only known examples are the polynomials $f(x)=x^2-a^2\in \bbF_q[x]$ with $q$ odd, $a\in\bbF_q$ of even order $m>4$ (where $\la a\ra$ is indeed a non-standard $f$-subgroup \cite[Proposition 2.3]{bn1}). In a subsequent paper, we will describe some new examples of degree 2.
\section{Acknowledgments}
We wish to thank Owen Brison and Eurico Nogueira for their careful reading of a draft version of this work. Their comments greatly helped to improve the paper. Part of this work has been initiated during a visit of the first author to Brison and Nogueira in 2013.
Both authors wish to acknowledge the support of Nanyang Technological University, Singapore, where the first
part of this research was carried out.
The research of H.D.L. Hollmann is in part supported by the Singapore National Research
Foundation under Research Grant NRF-CRP2-2007-03, 
and in part by the Estonian Research Council grant PRG49.
\end{document}